\documentclass[11pt,notitlepage,twoside,a4paper]{amsart}
\usepackage{amsmath}
\usepackage{amsfonts}
\usepackage{amssymb}
\usepackage{amscd}
\usepackage{amsthm}
\usepackage{amsbsy}
\usepackage[french, british]{babel}

\usepackage{amsmath,amssymb,enumerate}

\usepackage{epsfig,fancyhdr,color}

\usepackage{amssymb}
\usepackage{amsmath,amsthm}
\usepackage{latexsym}
\usepackage{amscd}
\usepackage{psfrag}
\usepackage{graphicx}
\usepackage[latin1]{inputenc}
\usepackage[all]{xy}
\usepackage[mathcal]{eucal}

\definecolor{NoteColor}{rgb}{1,0,0}

\renewcommand{\textsc}{\textcolor{red}}

\newtheorem*{theorem 1}{\rm\bf Proposition 1}
\newtheorem*{theorem 2}{\rm\bf Proposition 2}

\theoremstyle{definition}

\theoremstyle{remark}

\def\interieur#1{\mathord{\mathop{\kern 0pt #1}\limits^\circ}}

\begin{document}
 
\title[Map drawing and foliations of the sphere]{Map drawing and foliations of the sphere}

\author{Athanase Papadopoulos}
\address{Athanase Papadopoulos,  Universit{\'e} de Strasbourg and CNRS,
7 rue Ren\'e Descartes,
 67084 Strasbourg Cedex, France}
\email{papadop@math.unistra.fr}

\date{\today}
\maketitle

\begin{abstract}

 I will consider some questions related to Euler's work on cartography and its consequences, 
 in which the foliations of the sphere by meridians and parallels play important roles. 

\bigskip



\end{abstract}

 
%
%
%
%
%
%

   \section{Introduction}
    I will consider some questions related to Euler's work on the drawing of geographical maps, 
 in which the foliations of the sphere by meridians and parallels play an important role. 
This article can be regarded as a sequel to my article \cite{2016-Tchebyshev} in which I talked about the works of Euler and Chebyshev on  geography. It is also the occasion to straighten out a statement often made in the literature  concerning Euler's contribution to cartography.

\section{On Euler's ``perfect" mappings}
 I will start with a few remarks on Euler's memoir  \emph{De repraesentatione superficiei sphaericae super plano} (On the representation of  spherical surfaces on a plane) \cite{Euler-rep-1777}, presented to the Saint Petersburg Academy of Sciences on September 4, 1775, and published in the 1777 volume of the \emph{Acta Academiae Scientarum Imperialis Petropolitinae}. 
  
  This memoir is very poorly quoted in the literature. On the one hand, it is referred to for a result which is attributed to Euler, although this result was known since Greek antiquity, and its proof follows immediately from results in spherical geometry that were known in that period. On the other hand, there are interesting results that Euler proved in this memoir which, to my knowledge, are never mentioned in the papers or books on mathematical cartography. 
  
 The result that is usually attributed to Euler in relation with the paper \cite{Euler-rep-1777} says that there is no ``perfect" map from a subset of the sphere to the Euclidean plane. A confusion is entertained by the fact that the  word ``perfect" was used by Euler in his paper without a proper definition. In the (relatively large number of) papers in which Euler's memoir is quoted, it is said that Euler proved  in this memoir that there is no map from a subset of the sphere onto the Euclidean plane that preserves distances up to scale. The reason for this situation is that, as it often happens in ``historical" papers, authors get their information from papers written by other authors on the subject without bothering to look into the original sources and try to understand them.
 Let us give a few examples of such quotes.

    T. Feeman, in his book \emph{Portraits of the Earth}, discusses the existence of a map from a portion of the sphere onto the plane which has a fixed scale, that is, as the author puts it, a map such that the the quantity
   \[\frac{\text{distance between two points on the globe}}{\text{distance between their images}}
   \] is constant, i.e. independent of the chosen pair of points. The author writes  \cite[p. 25]{Feeman} (2002):
   ``Over the years various attempts by cartographers to solve this problem resulted in some ingenious, if flawed, maps. Finally, in 1775, Leonhard Euler (1707--1783), the leading mathematician of his day and one of the most important mathematical figures of all time, presented to the St. Petersburg Academy of Sciences a paper entitled \emph{On representations of a spherical surface on the plane} in which he proved conclusively that such a map could not exist." 
   
      In the article \emph{Dallo spazio come contenitore allo spazio come rete}, C. Corrales Rodriganez writes \cite[p. 125]{Rodriganez}:  ``The mathematician Leonhard Euler, in his article \emph{De representatione superficiei sphaericae super plano}, published in the eighteenth century, proved that no part of the Earth can be reproduced over a plane surface without deformation. Euler's theorem says that the perfect map does not exist."\footnote{[Il matematico Leonard Euler] nel suo articolo \emph{De representatione superficiei sphaericae super plano}, pubblicato nel Settecento, ha provato che nessuna parte della Terra può essere riprodotta su une superficie piana senza deformazione. Il teorema di Eulero dice che la carta perfetta non esiste.} Incidentally, the author adds that ``Euler's theorem has pushed the mathematical cartographers to study spherical geometry and trigonometry as a subject in itself, independent of Euclidean geometry."\footnote{Il teorema di Eulero ha spunto i cartografi matematici a studiare la geometria sferica e la trigonometria come materie a sé, indimendenti dalla geometria euclidea.} It is not clear to what mathematical cartographers the author is referring to, but Euler himself was a cartographer, and he had published several works on spherical geometry and spherical trigonometry several decades before he wrote this memoir.\footnote{An edition of Euler's and his collaborators on spherical geometry will appear in the book \cite{Caddeo2}.} Furthermore, long before Euler, Ptolemy (2nd c.A.D.), one of the most famous mathematical cartographers of all times, was also thoroughly involved in spherical geometry and spherical trigonometry. As a matter of fact, the field of spherical geometry remains poorly known to historians of mathematics. In the book \emph{Portraits of the Earth}  which we already mentioned, the author claims (p. 25) that spherical geometry can be developed axiomatically as a Euclidean geometry in which Euclid's fifth postulate is replaced by a postulate saying that any two lines intersect (in two points). This statement is not correct. The confusion is probably caused by the fact that the axioms of hyperbolic geometry (which, together with Euclidean and spherical geometry forms the three ``classical" geometries, or the geometries of constant curvature) are precisely those of Euclidean geometry with the fifth postulate replaced by its negation. Spherical geometry may be developed axiomatically, but such a set of axioms cannot be so simply obtained from those of Euclidean geometry.
      
      Let us continue with our citations of Euler's result.
      
   R. Osserman, in his paper 
 \emph{Mathematical mapping from Mercator to the millennium} \cite[p. 234]{Osserman} (2004) attributes the following theorem to Euler: \emph{It is impossible to make an exact scale map of any part of a spherical surface}. By an ``exact scale map", Osserman means a map that preserves distances up to scale, and he refers again to Euler's paper \cite{Euler-rep-1777}.  

P. Robinson, in a paper titled \emph{The sphere is not flat} \cite{Robinson} (2006), writes the following:
``The theorem of our title asserts that there is no isometric (that is, distance-preserving) function from the sphere (or indeed from any of its nonempty open subsets) to the Euclidean plane; more generally, there is no isometry to any Euclidean space. This theorem may be traced back to Euler, in his \emph{De repraesentatione superficiel sphaericae super plano} of 1778." 
 
 The same poor attribution to Euler occurs in the chapter titled \emph{Curvature and the notion of space}  in the book \emph{Mathematical Masterpieces} (2007) by  A. Knoebel, J.  Lodder, R.  Laubenbacher and D. Pengelley \cite[p. 163]{KLLP}, where the authors write:
 ``In the paper presented to the St. Petersburg Academy of Science in 1775, \emph{De repraesentatione superficiei sphaericae super plano}, Euler proved what cartographers had long suspected, namely, the impossibility of constructing a flat map of the round world so that all distances on the globe are proportional (by the same constant of proportionality) to the corresponding distances on the map."

   J.  Gray, in his book  \emph{Simply Riemann} \cite{Gray} which appeared in 2020, says the following, about Euler's work on cartography:  
   ``Euler used all his analysis to prove that every cartographer suspected: that there could be no map of the Earth's surface onto a plane that is accurate in every respect. Some maps send curves of shortest length on the sphere to straight lines in the plane; there are maps that send equal angles to equal angles, and there are maps that scale all areas by the same amount. But there can be no map that does all of these at once."

  Naturally, popular science authors get their information from mathematicians' writings when they understand them: In the Spanish daily newspaper \emph{La Vanguardia}, on 26 March 2017, in an article titled \emph{Un mundo, tres mapas}, the author, A. Molins Renter, writes: ``Passing from a spherical geometric form to a plane support, two-dimensional and usually of rectangular shape, results in the fact that something is always lost in the translation, as the Swiss mathematician and physicist Leonhard Euler already demonstrated in 1778, in his work \emph{De repraesentatione superficiei sphaericae super plano}."\footnote{Pasar de una forma geométrica esférica a un soporte plano, bidimensional y normalmente con forma rectangular provoca que algo se pierda siempre en la translaci\'on, como ya demostr\'o el matem\'atico Leonhard Euler en 1778, en su obra \emph{De  repraesentatione superficiei sphaericae super plano}.}
  
  One could give are many other examples.

  In fact, the statement attributed by all these authors to Euler was obviously known to him, but it does not convey the slightest idea of the results obtained in the memoir quoted, which are much stronger and much more interesting than what all these authors claim. Furthermore, as I said, the result that they quote was known since the 1st-2nd century A.D., since it follows as a corollary from several results contained in Menelaus' \emph{Spherics} on the geometry of spherical triangles. For instance, it is an immediate consequence of the result saying that the angle sum in a spherical triangle is always greater than 2 right angles (this is Proposition 12 in \cite{RR2}), or from the comparison result saying that in any spherical triangle $ABC$, if $D$ and $E$ denote the
midpoints of $AB$ and $B$ respectively and if $DE$ is the shortest arc joining them. Then $DE > AC/2$ (Proposition 27 in \cite{RR2}). A local isometry between an open region of the sphere and an open region of the Euclidean plane would preserve the two properties of triangles, and obviously neither of them is satisfied by a Euclidean triangle. 
  
    At the end of Section 9 of his paper, Euler writes: ``It is proved through computation that a perfect mapping of the Sphere
onto the plane is not possible." But since he did not give  the definition of a \emph{perfect map}, the meaning of this sentence should be understood in context, that is, by following the arguments that lead to it. 

Recently, C. Charitos and I. Papadoperakis wrote a paper titled \emph{On the non-existence of a perfect map
from the 2-sphere to the Euclidean plane} \cite{CP} in which they give a precise statement of Euler's result and provide a detailed proof  of it.  

To state Euler's result correctly, we call a
 map $f$ from a region  $S$ of the 2-sphere to the Euclidean plane \emph{perfect} if every point in the domain has a neighborhood on which the following two conditions hold:
  
  \begin{enumerate}
  \item $f$ sends meridians and parallels to two fields of lines that make mutually the same angles;
  \item $f$ preserves distances infinitesimally along the meridians and the parallels.
  \end{enumerate}
Thus, a perfect map sends the meridians and parallels to two line fields that are orthogonal.
Furthermore, a perfect map preserves globally the length element along the meridians and the parallels.    One should note here that on the spherical globe, the meridians are geodesics but the parallels are not. The fact that distances are preserved infinitesimally along the meridians implies immediately that the distances between points on these lines are preserved. It also follows, although not so immediately, that distances between points on parallels is preserved by a perfect map.
  
  The idea of Euler's proof was to translate these geometrical conditions into a
system of partial differential equations and to show that this system has no solution.
  
    Furthermore, Euler, in his paper, after showing the non-existence of a perfect map, proves several other results. He declares that since perfect maps do not exist, one has to look for best approximations.  He writes:  ``We are led to consider representations which are not similar, so that the spherical figure differs in some manner from its image in the plane." He then examines several particular projections of the sphere, searching systematically for the partial differential equations that they satisfy.  He considers several classes of maps: conformal maps (which he calls ``similitudes on the small scale"), area-preserving maps, and maps where the images of all the meridians are perpendicular to a given axis while those of all parallels are parallel to it.  He gives examples of maps satisfying each of the above three properties and in each case he studies their distance and angle distortion.
     
   \section{On the action of a geographical map on the foliations by parallels and meridians}  An important feature of Euler's memoir \cite{Euler-rep-1777} and of the other memoirs that we shall consider below is that most of the properties of the geographical maps that are requested are formulated in terms of how these maps transform the two geographically most famous foliations of the sphere, namely, the foliations by parallels and  by meridians.
   
   Let us recall that when the surface of the Earth is considered to be a sphere, the \emph{parallels} are the family of circles that are equidistant from the equator. The latter is the great circle that is perpendicular to the rotation axis of the Earth, the one that separates the Northern hemisphere from Southern one. In the geometry of the sphere, the parallels are geometric circles, that is, equidistant points from a center, which is either the North or the South pole (a circle on the sphere has two centers). Furthermore, the parallels are small circles, that is, intersections of the sphere with planes that do not pass through the center, except the equator itself, which is also considered as a parallel (at zero distance from itself) and which is a great circle, the intersection of the sphere with a plane passing through the center. The difference is important because on the sphere great circles are geodesics whereas small circles are not. This foliations has two singular points, situated at the North and South poles. 
      
The second foliation that is used in the paper  \cite{Euler-rep-1777} is the foliation by \emph{meridians},  whose leaves are the great circles perpendicular to the equator, or, equivalently, the great circles that pass through the North and South poles. Unlike the parallels, the meridians are all great circles, and therefore geodesics. It has two singular points, situated at the North and South poles. 

\begin{figure}[htbp]
\centering
\includegraphics[width=7cm]{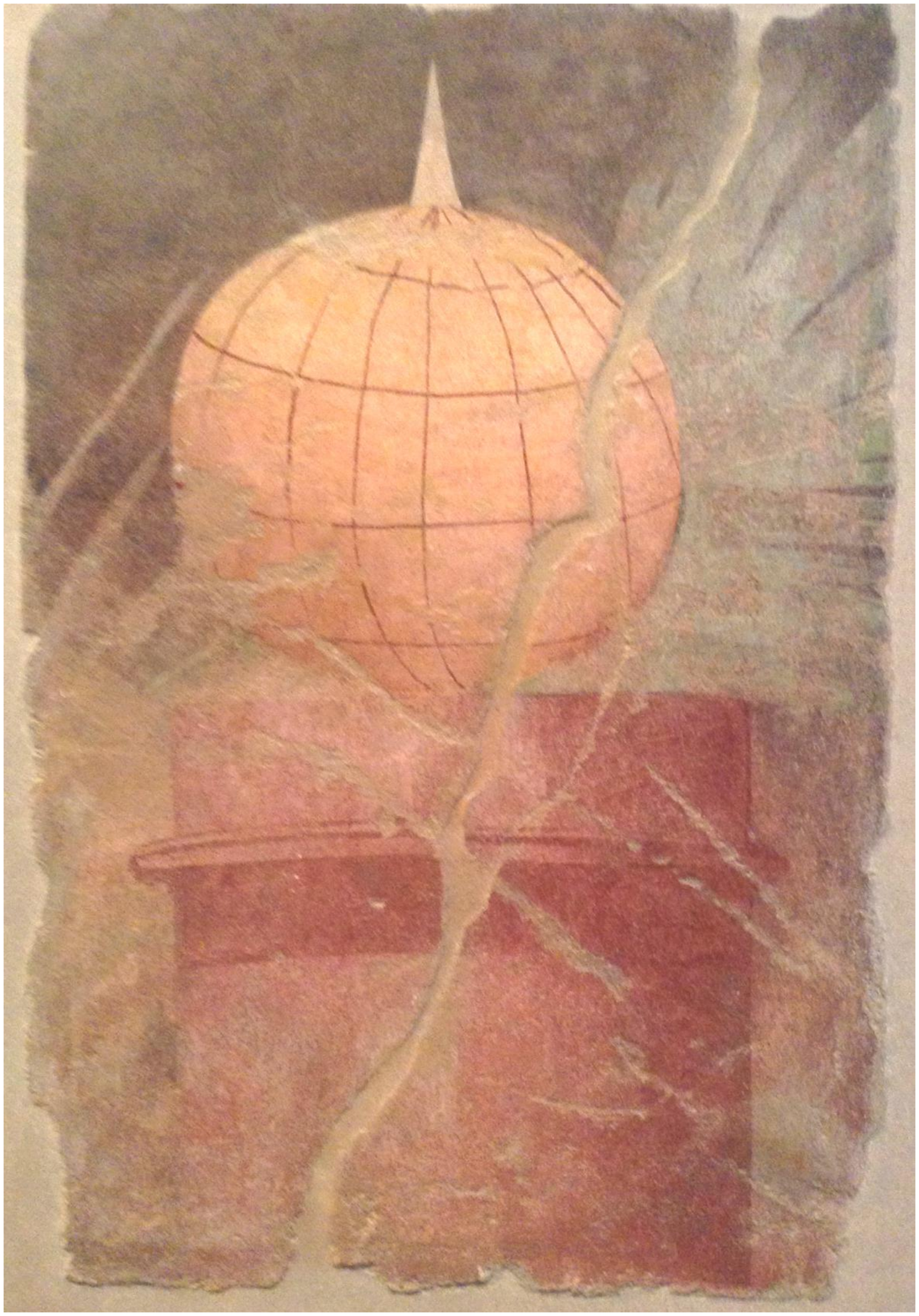}
\caption{\small A celestial globe with its two perpendicular foliations by parallels and meridians. Wall painting fragment from the first century A.D.  Metropolitan Museum of Art, New York,  Department of Greek and Roman Art. (Photo A. Papadopoulos.)} \label{Met}
\end{figure}

Since the time of ancient Greek geography, the foliations by parallels and meridians play an important role in map drawing, for representing regions of the Earth, but also for maps of the celestial sphere.
 Figure 
    \ref{Met} is a reproduction of a representation of a celestial globe dating from the 1st century A. D. on which the foliations by parallels and meridians are drawn. The same picture could serve for the representation of the Earth.
           
      The properties of the images of the meridians and parallels are important factors in several known projections. For instance, under a stereographic projection centered at the North pole, the parallels are sent to concentric circles centered at the image of the South pole, while the meridians are sent to straight lines meeting at the North pole. The projection from the center of the sphere to a plane tangent to the South pole,  known as the gnomonic projection, and which was used since the times of Thales, has a similar property: parallels are sent to circles centered at the South pole and meridians are sent to straight lines  passing through this pole. 
      
\begin{figure}
\centering
 \includegraphics[width=\linewidth]{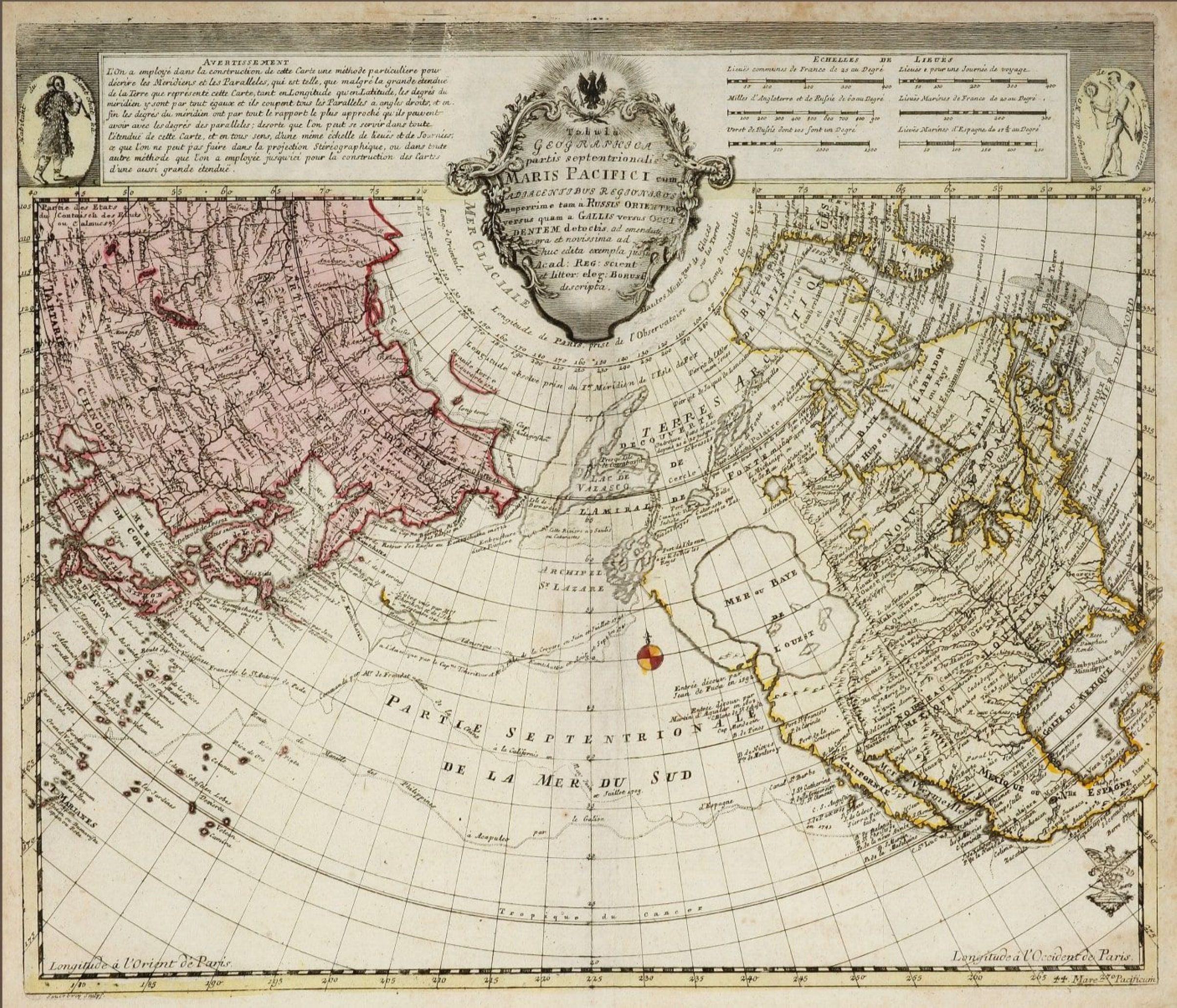}    \caption{\small A map of the Northern Pacific, with the Eastern part of Asia and the Northern part of America, from Euler's \emph{Atlas Geographicus} (Berlin, 1753)}   \label{fig:Euler-Map41} 
\end{figure}

      The foliations by parallels and by meridians are perpendicular. Their images  by a geographical map are usually drawn on the map, see e.g. Figures \ref{fig:Euler-Map41} and \ref{Delisle-1745}.

Two other memoirs by Euler on geography were published the same year as his memoir \emph{De repraesentatione superficiei sphaericae super plano}, namely,   \emph{De proiectione geographica superficiei sphaericae} \cite{Euler-pro-1777} and \emph{De proiectione geographica Deslisliana in mappa generali imperii russici usitata} \cite{Euler-pro-Desli-1777}. I would like to make a few comments on the latter.

 The memoir \cite{Euler-pro-Desli-1777} is concerned with a projection from the sphere that was used by Joseph-Nicolas Delisle who, during several years, was the main geographer and the director of the astronomical department of the Saint Petersburg  Academy of Sciences. He was in charge of drawing new and precise maps of the Russian Empire. From 1735 to 1740, Euler assisted Delisle in this work until he became himself the head of the geography department of the Academy, after  a conflict emerged between 
Delisle and the Academy's administration, in relation with the so-called \emph{Atlas Russicus} (the ``Russian Atlas"), a project initiated by Peter the Great and of which Delisle was in charge, and which he kept postponing. In 1740, the responsibility of the Russian Atlas was taken away from him and given to Euler. 
  
 In his memoir, Euler starts by reviewing the main properties of a stereographic projection used by the geographer Johann Matthias Hasius. The latter had published in Nuremberg, in 1739, a map  of Russia known under the name ``Imperii Russici et Tatariae universae tam majoris et asiaticae quam minoris et europaeae tabula" (Geographical map of the Russian Empire and of Tataria, both large and small, in Europe and Asia).
 Euler mentions properties of the images of the two foliations by parallels and meridians, in particular that these images intersect at right angles (in fact, the map is conformal, that is, it is angle-preserving). Euler then reviews the inconveniences of this projection: length is highly distorted in the large, especially for maps that represent large regions of the Earth and the images of the meridians are not evenly curved on the geographical map, even though these lines are circles. In particular, the province of Kamchatka is distorted by a factor of four, compared to another region at the center of the map.

In \S 5, of his memoir, Euler states the following  four properties that are required from an ideal geographical map:
 (i) the images of the meridians are straight lines; (ii)
the  degrees of latitudes do not change along meridians;
(iii) the images pf the parallels meet the images of the meridians at right angles; 
(iv) at each point of the map, the ratio of the degree on the parallel to the degree on the meridian is the same as on the sphere.
He then declares that since this cannot be achieved one may request,  instead of the last condition,  that the deviation of the degree of latitude to the degree of longitude at each point from the true ratio be
 as small as possible (ideally, this error should not be noticeable).

He then recalls the construction of Delisle's map.

\begin{figure}
\centering
 \includegraphics[width=\linewidth]{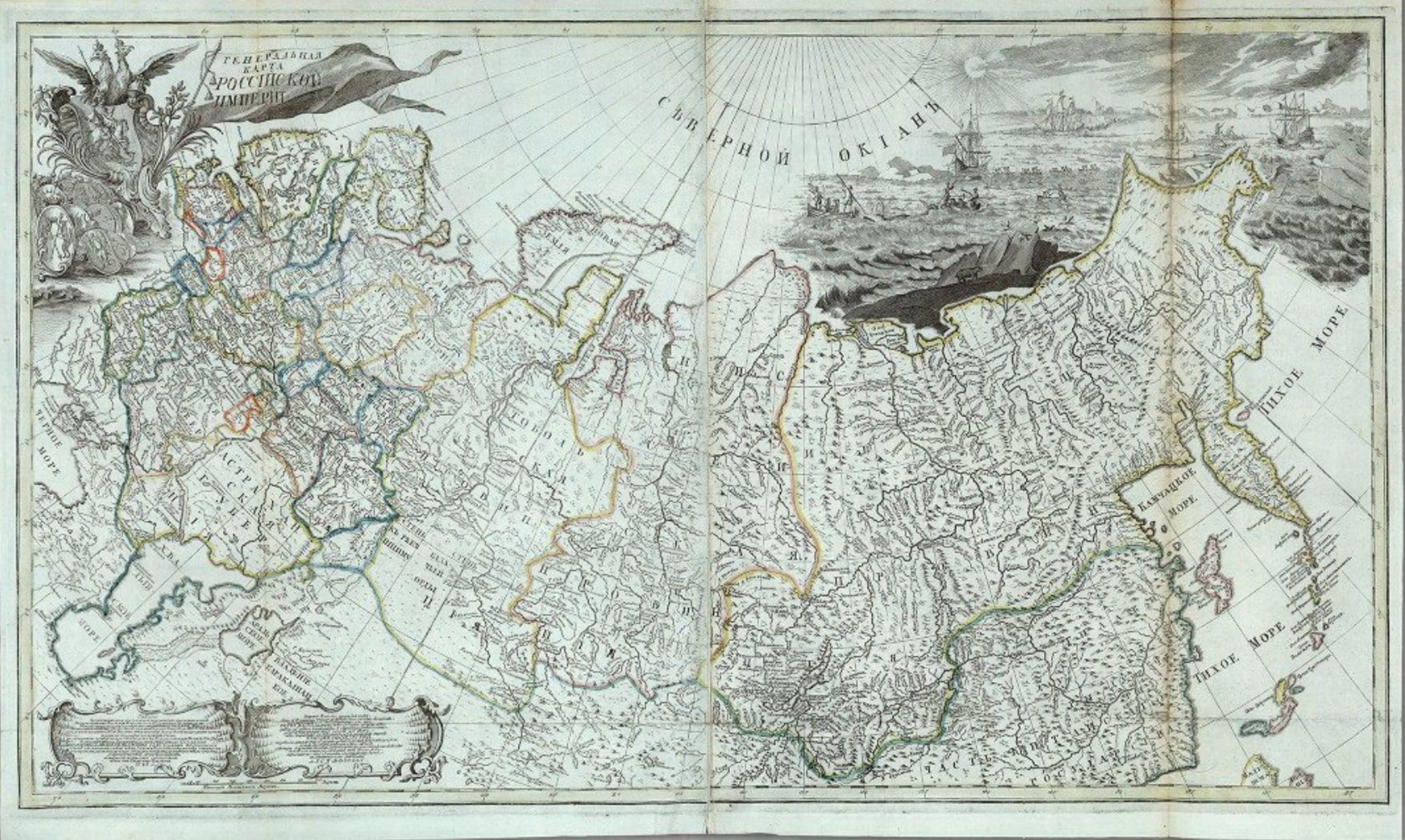}    \caption{\small A map of the Russian Empire by Joseph-Nicolas Delisle (Saint Petersburg, 1745)}   \label{fig:Delisle}
\end{figure}

In this map, one first chooses  two outer parallels that contain the region that is to be represented. In the case of the Russian Empire, these outermost parallels are chosen to be those at $40^{\mathrm{o}}$ and $70^{\mathrm{o}}$ of altitude.  Then,  to chooses two other special parallels on which the ratios of the degrees of latitude to the degrees of longitude will be represented by their exact values.  Euler writes that the question becomes that of choosing these two new parallels in such a way that the maximum deviation of the  ratios of the degrees of latitude and longitude
over the entire map is minimized. He writes that Delisle found that the optimal choice of these parallels is to take them equidistant from the central parallel of the map and from the  outermost parallels chosen.
Besides, distances should be preserved on all meridians, and the maximum of a certain deviation over the entire map must be minimized.   
Figure \ref{fig:Delisle} reproduces a map drawn by Delisle using his method.

Starting in \S 7,  Euler presents a mathematical construction of a family of  straight lines representing meridians which are at distance one degree from each other. He
notes that one advantage of Delisle's projection is that while meridians are represented by straight lines, the images of the other great circles do not deviate considerably from straight lines (\S 22 of \cite{Euler-pro-Desli-1777}), and he gives some precise estimates of this  deviation (\S 23ff).

In \S 10--16 of his memoir, Euler gives the mathematical details that show the difference between this representation and the real situation at the extreme points he started with, and in \S 17--23, he makes the actual computations in the  special case where the map is that of the Russian Empire. In the last sections of his memoir,  he studies the images of great circles on the sphere and he shows that the difference between these images and circular arcs is not noticeable. He computes the radius of such an arc, which he finds to be very large and he concludes that the shortest lines on the map do not differ sensibly from straight lines.

 In conclusion, let us stress once again that in Euler's treatment of Delisle's projection, like in the theorem we stated above on perfect maps, the important requirements concern the behavior of the images of the two  foliations defined by the parallels and the meridians.

  Besides the memoir  \cite{Euler-pro-Desli-1777} in which he described Delisle's method of drawing geographical maps, Euler explained the same method in the \emph{Atlas Geographicus omnes orbis terrarum regiones in XLI tabulis exhibens} \cite{Atlas} that was published by the \emph{Académie Royale des Sciences et Belles Lettres de Prusse}, in the year 1753, in Berlin, where he worked for 25 years, between his two stays in Saint Petersburg. The atlas contains 45 maps, and was edited under the direction of Euler who also wrote its preface, which is dated May 13, 1753.  Several projections are used in this atlas, which is concerned only with large parts of the Earth. In all these projections, the meridians are perpendicular to the parallels.

The map in Figure \ref{fig:Euler-Map41} Figure \ref{fig:Euler-Map41} is extracted from this atlas. It is the last one in the series, and it is drawn using Delisle's method. Euler, in the preface, comments on this method. Euler writes that Delisle's method  seems to him the most appropriate for a proper representation of these Northern regions of the terrestrial globe. 
He recalls that in this representation, the meridians are straight lines and all their degrees are equal: the images of two meridians that are distant apart by one degree converge in such a manner that at two altitudes that are chosen in advance, the ratio of the degrees of longitude to the degrees of latitude are the same ratio as in reality.  This is the property that he presents in his memoir  \cite{Euler-pro-Desli-1777} that we noted above. For the Russian Empire,  after the choice of the two outer parallels  at $40^{\mathrm{o}}$ and $70^{\mathrm{o}}$, the two parallels  which are at the same elevation from the extremities of the region that is represented as well as from its center are those 
at $47^{\mathrm{o}}  30'$ and t $62^{\mathrm{o}}  30'$.
  Under these two altitudes, the ratios between the degrees of longitude and latitude are  accurate on the map. At the other locations, they are almost accurate (the difference is not noticeable). Besides, in this representation, all the meridians (which are straight lines) merge at a point, although this point is not the North pole; it is at a distance which would correspond to 7 degrees farther than this pole. From this point as center, the images of the parallels are circles. Euler writes that one should not regard as a shortage of this map the fact that the center in which all the meridians intersect is so far from the pole, nor the fact that on this map the parallels, which form semi-circles, do not occupy 180$^{\mathrm{o}}$ in longitude, but much more, sometimes even up to 250$^{\mathrm{o}}$.

    Lagrange, whose name is associated to the one Euler in several respects, already stressed in his paper  on the construction of geographical maps the fact that the only thing we have to do  drawing a geographical map is to specify the images of meridians and parallels according to a certain rule (see  \cite[p. 640]{Lagrange1779}). This simple remark was at the basis of the the development of modern mathematical cartography. We have tried to convey this idea by mentioning examples from the works of Euler and Delisle, but others may be found in works of Lambert, Gauss, Bonnet and others. The forthcoming book \cite{Caddeo2} contains a section on cartography at the epoch of Euler.
 
\bigskip

   \noindent {\it Acknowledgements.}
   This paper is based on a talk I gave in Obninsk, on May 14, 2019, at a conference commemorating Pafnouti Chebyshev. I would like to thank the organizer, Valerii Galkin, for inviting me to that conference. I would also like to thank Alena Zhukova for her corrections on a preliminary version of this paper.

\end{document}